\documentclass[preprint,12pt,authoryear]{elsarticle}

\usepackage{amssymb}
\newtheorem{Lem}{Lemma}

\newtheorem{Theor}{Theorem}

\newcommand{\R}{\ensuremath{\mathbb{R}}}

\newcommand{\cqfd}{\hfill $\square$}
\usepackage[usenames,dvipsnames]{xcolor} 

\usepackage{setspace}
\journal{some journal}

\begin{document}

\begin{frontmatter}

%% Title, authors and addresses

%% use the tnoteref command within \title for footnotes;
%% use the tnotetext command for theassociated footnote;
%% use the fnref command within \author or \address for footnotes;
%% use the fntext command for theassociated footnote;
%% use the corref command within \author for corresponding author footnotes;
%% use the cortext command for theassociated footnote;
%% use the ead command for the email address,
%% and the form \ead[url] for the home page:
%% \title{Title\tnoteref{label1}}
%% \tnotetext[label1]{}
%% \author{Name\corref{cor1}\fnref{label2}}
%% \ead{email address}
%% \ead[url]{home page}
%% \fntext[label2]{}
%% \cortext[cor1]{}
%% \address{Address\fnref{label3}}
%% \fntext[label3]{}

\title{On the Maximal Halfspace Depth of Permutation-invariant Distributions on the Simplex} %Majorization}

%% use optional labels to link authors explicitly to addresses:
%% \author[label1,label2]{}
%% \address[label1]{} 
%% \address[label2]{}
 
\author{Davy Paindaveine\fnref{DP2}}

\author{Germain Van Bever}  

\fntext[DP2]{Corresponding author (dpaindav@ulb.ac.be)}

\address{Universit\'{e} libre de Bruxelles, ECARES and D\'{e}partement de Math\'{e}matique, CP 114/04, 50, Avenue F.D. Roosevelt, B-1050 Brussels, Belgium}

\begin{abstract}
We compute the maximal halfspace depth for a class of permutation-invariant distributions on the probability simplex. The derivations are based on stochastic ordering results that so far were only showed to be relevant for the Behrens-Fisher problem. 
\end{abstract}

\begin{keyword}
$\alpha$-unimodality \sep Dirichlet distribution \sep Halfspace depth \sep Majorization \sep Stochastic ordering
%% keywords here, in the form: keyword \sep keyword

%% PACS codes here, in the form: \PACS code \sep code

%% MSC codes here, in the form: \MSC code \sep code
%% or \MSC[2008] code \sep code (2000 is the default)

\end{keyword}

\end{frontmatter}

%% \linenumbers

%% main text
\section{Introduction}
\label{SecIntro}
 
Denoting as~$\mathcal{S}^{k-1}:=\{x\in\R^k : \|x\|^2=x'x=1\}$ the unit sphere in~$\R^k$, the \cite{Tuk1975} halfspace depth
$
H\!D(\theta,P)
 = 
\inf_{u\in\mathcal{S}^{k-1}}  
P\big[ u'(X-\theta) \geq 0 \big]
$
measures the centrality of the $k$-vector~$\theta$ with respect to a probability measure~$P=P^X$ over~$\R^k$. Any probability measure~$P$ admits a deepest point, that generalizes to the multivariate setup the univariate median; see, e.g., Proposition~7 in \cite{RouRut1999}. 
Parallel to the univariate median, this deepest point is not unique in general. Whenever a unique representative of the collection~$C_P$ of $P$'s deepest points is needed, the \emph{Tukey median}~$\theta_P$, that is defined as the barycentre of~$C_P$, is often considered. The convexity of~$C_P$ (see, e.g., the corollary of Proposition~1 in \citealp{RouRut1999}) implies that~$\theta_P$ has maximal depth. The depth of~$\theta_P$ is larger than or equal to~$1/(k+1)$; see Lemma~6.3 in \cite{DonGas1992}.     

In this paper, we determine the Tukey median and the corresponding maximal depth for a class of permutation-invariant distributions on the probability simplex~$\mathcal{D}_k:=\{x\in\R^k: x_1,\ldots,x_k \geq 0, \ \sum_{\ell=1}^k x_\ell = 1\}$. The results identify the most central location for some of the most successful models used for compositional data. They have also recently proved useful in the context of depth for shape matrices; see \cite{PVB17b}. 

%The derivations are based on stochastic ordering results that, to the best of our knowledge, were only used so far in the framework of the Behrens-Fisher problem; see \cite{Haj1962}, \cite{Law1968}, and \cite{EatOls72}.
%We are not aware of any other method to compute the maximum halfspace depth for the distributions considered. 

The outline of the paper is as follows. In Section~\ref{SecDef}, we define the class of distributions we will consider and state a stochastic ordering result on which our derivations will be based. In Section~\ref{SecMain}, we state and prove the main results of the paper. Finally, in Section~\ref{SecSimu}, we illustrate the results through  numerical exercises and we shortly comment on open research questions.

%%%%%

\section{Preliminaries}
\label{SecDef}
 
Let~$\mathcal{F}$ be the collection of cumulative distribution functions~$F$ such that (i)~$F(0)=0$ and~(ii)~$F$ is concave on~$(0,+\infty)$. In other words,~$\mathcal{F}$ collects the cumulative distribution functions of random variables that are (i) almost surely positive and (ii) unimodal at~$0$. Any~$F$ in~$\mathcal{F}$ admits a probability density function~$f$ that is non-increasing on~$(0,+\infty)$. 

For an integer~$k\geq 2$ and~$F$ in~$\mathcal{F}$, consider the random $k$-vector 
\begin{equation}
\label{defdir}	
X_k
=
(X_{k1},\ldots,X_{kk})'
=
(V_{k1},\ldots,V_{kk})'/\,{\textstyle{\sum_{\ell=1}^k}} V_{k\ell}
,
\end{equation}
where~$V_{k1},\ldots,V_{kk}$ are mutually independent and have cumulative distribution function~$F$. The corresponding probability measure over~$\R^k$ will be denoted as~$P_{k,F}$. Obviously, the random vector~$X_k$ takes its values in the probability simplex~$\mathcal{D}_k$.
This includes, for example, the Dirichlet distribution with parameter~$(\alpha,\ldots,\alpha)'\in\R^k$, obtained for the cumulative distribution function~$F=F_\alpha$ of the Gamma$(\alpha,\frac{1}{2})$ distribution (that corresponds to the probability density function~$x\mapsto x^{\alpha-1}\exp(-x/2)/(2^\alpha\Gamma(\alpha))$ on~$(0,+\infty)$, where $\Gamma$ is the Euler Gamma function). The unimodality constraint in~(ii) above imposes to restrict to~$\alpha\leq 1$. Note that, irrespective of~$F$, the mean vector of~$X_k$ is~$\mu_k=k^{-1}1_k=k^{-1}(1,\ldots,1)'$.

To state the stochastic ordering result used in the sequel, we need to introduce the following notation. For $n$-vectors~$a,b$ with~$\sum_{\ell=1}^n a_\ell=\sum_{\ell=1}^n b_\ell$, we will say that \emph{$a$ is majorized by~$b$} if and only if, after permuting the components of these vectors in such a way that~$a_1\geq a_2\geq \ldots\geq a_n$ and~$b_1\geq b_2\geq \ldots\geq b_n$ (possible ties are unimportant below),  
$
\sum_{\ell=1}^r a_\ell 
\leq
\sum_{\ell=1}^r b_\ell 
$,
for any~$r=1,\ldots,n-1$;
see, e.g., \cite{Maretal2011}. For random variables~$Y_1$ and~$Y_2$,
%, with respective cumulative distribution functions~$F_1$ and~$F_2$, 
we will say that~$Y_1$ is stochastically smaller than~$Y_2$ ($Y_1\leq_{\rm st}Y_2$) if and only if~$P[Y_1>t]\leq P[Y_2>t]$ for any~$t\in\R$. To the best of our knowledge, the following stochastic ordering result so far was only used in the framework of the Behrens-Fisher problem; see \cite{Haj1962}, \cite{Law1968}, and \cite{EatOls72}.

 \begin{Lem}[\citealp{EatOls72}]
 	\label{Olshen}
Let~$W$ be a random variable with a cumulative distribution function in~$\mathcal{F}$. Let~$Q_1,\ldots,Q_n$ be exchangeable positive random variables that are independent of~$W$. Then, 
%$$
%\frac{
%W 
%} 
%{
%\sum_{\ell=1}^{n} a_\ell Q_\ell
%}
%\leq_{\rm st} 
%\frac{
%W
%}
%{
%\sum_{\ell=1}^{n} b_\ell Q_\ell
%}
%$$ 
$
W 
/
(\sum_{\ell=1}^{n} a_\ell Q_\ell)
\leq_{\rm st} 
W
/
(\sum_{\ell=1}^{n} b_\ell Q_\ell)
$ 
for any~$a,b\in\R^n$ such that~$a$ is majorized by~$b$.  
  \end{Lem}
  
In \cite{EatOls72}, the result is stated in a vectorial context that requires the $\alpha$-unimodality concept from \cite{OlsSav1970}. In the present scalar case, the minimal unimodality assumption in \cite{EatOls72} is that~$\sqrt{W}$ is $2$-unimodal about zero, which, in view of Lemma~2 in  \cite{OlsSav1970}, is strictly equivalent to requiring that~$W$ is unimodal about zero.

%%%%%

\section{Main results}
\label{SecMain}

%${\rm Beta}(1/2,(k-1)/2)$. 

Our main goal is to determine the Tukey median of the probability measure~$P_{k,F}$ and the corresponding maximal depth. Permutation invariance of~$P_{k,F}$ and affine invariance of halfspace depth allows to obtain %the following result. 

\begin{Theor}
\label{thetheoremmax}
The Tukey median of~$P_{k,F}$ is~$\mu_k=k^{-1}1_k$. 
%	$H\!D(\mu_k ,P_k)\geq 	H\!D(\theta ,P_k)$ for any~$\theta\in \mathcal{D}_k$.
	\end{Theor}

{\sc Proof of Theorem~\ref{thetheoremmax}}. 
Let~$\theta_*=(\theta_{*1},\ldots,\theta_{*k})'$ be a point maximizing~$H\!D(\theta ,P_{k,F})$ and let~$\alpha_*=H\!D(\theta_*,P_{k,F})$ be the corresponding maximal depth. Of course, $\theta_*\in\mathcal{D}_k$ (if~$\theta\notin \mathcal{D}_k$, then~$H\!D(\theta ,P_{k,F})=0$). Denote by~$\pi_i$, $i=1,\ldots,k!$, the $k!$ permutation matrices on $k$-vectors. By affine invariance of halfspace depth and permutation invariance of~$P_{k,F}$, all~$\pi_i\theta_*$'s have maximal depth~$\alpha_*$ with respect to~$P_{k,F}$. Now, for any~$\ell=1,\ldots,k$, 
$$
\frac{1}{k!} \sum_{i=1}^{k!} 
(\pi_i\theta_*)_\ell
=
\frac{1}{k!} \sum_{\ell=1}^{k} 
(k-1)! \theta_{*\ell}
=
\frac{1}{k} \sum_{\ell=1}^{k}  \theta_{*\ell}
=
\frac{1}{k}
=
(\mu_k)_\ell.
$$
Since this holds for any~$\theta_*$ maximizing~$H\!D(\theta ,P_{k,F})$, the result is proved. 
%Quasi-concavity of halfspace depth then yields
%\begin{eqnarray*}
%H\!D(\mu_k,P_{k,F})
%&\!\!=\!\!&
%H\!D\bigg( \frac{1}{k!} \sum_{i=1}^{k!} \pi_i\theta_* ,P_{k,F}\bigg)
%\\[2mm]
%&\!\! \geq \!\!&
%\min\big(H\!D( \pi_1\theta_* ,P_{k,F}),\ldots,H\!D( \pi_{k!}\theta_* ,P_{k,F})\big)
%=
%\alpha_*
%.
%\end{eqnarray*}
%Since~$\alpha_*=\max_{\theta\in\R^k} H\!D(\theta,P_{k,F})$ by definition, we have that~$H\!D(\mu_k,P_{k,F})\leq \alpha_*$, hence that~$H\!D(\mu_k,P_{k,F}) = \alpha_*$.  
\cqfd
 \vspace{3mm}

%%%%

Note that the unimodality of~$F$ about zero is not used in the proof of Theorem~\ref{thetheoremmax}, so that the result also holds at~$F_\alpha$ with~$\alpha\geq 1$. In contrast, the proof of the following result, that derives the halfspace depth of the Tukey median of~$P_{k,F}$, requires unimodality.

\begin{Theor}
\label{thetheorem}
Let~$X_k=(X_{k1},\ldots,X_{kk})'$ have distribution~$P_{k,F}$ with~$F\in\mathcal{F}$. Then, 
$
	H\!D(\mu_k,P_{k,F})
	=
%	h_k
%	$, 
%	where we let~$h_k
%=
 	P\big[
X_{k1}
\geq 
1/k
 \big]
.
 $
% where~$Y_k\sim {\rm Beta}(1/2,(k-1)/2)$.
\end{Theor}

The proof requires the following preliminary result.

 \begin{Lem}
 	\label{thelemma}
For any positive integer~$k$, let $
h_{k,F}
=
 	P\big[
X_{k1}
\geq 
1/k
 \big]
 $,
where~$X_k=(X_{k1},\ldots,X_{kk})'$ has distribution~$P_{k,F}$ with~$F\in\mathcal{F}$.
 Then, the sequence~$(h_{k,F})$ is monotone non-increasing. 
  \end{Lem}

  {\sc Proof of Lemma~\ref{thelemma}}. 
  Since~$
k^{-1} 1_{k}$ is majorized by the $k$-vector~$
(k-1)^{-1} (1,\ldots,1,0)'$,  
Lemma~\ref{Olshen} readily provides
\begin{eqnarray*}
\lefteqn{
\hspace{-8mm} 
h_{k+1,F}
=
  P
  \Bigg[
  V_{k+1,1}
  \geq 
  \frac{1}{k+1} 
 \sum_{\ell=1}^{k+1} V_{k+1,\ell}
  \Bigg]
=
  P
  \Bigg[
  \frac{k}{k+1} 
  V_{k+1,1}
  \geq 
  \frac{1}{k+1} 
 \sum_{\ell=2}^{k+1} V_{k+1,\ell}
  \Bigg]
}
  \\[2mm]
  & & 
  \hspace{10mm} 
=
  P
  \Bigg[
\frac{V_{k+1,1}}{  \frac{1}{k} 
 \sum_{\ell=2}^{k+1} V_{k+1,\ell} 
}
  \geq 
1
  \Bigg]
\leq
  P
  \Bigg[
\frac{V_{k+1,1}}{  \frac{1}{k-1} 
 \sum_{\ell=2}^{k} V_{k+1,\ell}
}
  \geq 
1
  \Bigg]
  \\[2mm]
  & & 
  \hspace{10mm} 
  =
  P
  \Bigg[
  V_{k1}
  \geq 
  \frac{1}{k-1} 
 \sum_{\ell=2}^{k} V_{k\ell} 
  \Bigg]
  =
  P
  \Bigg[
  V_{k1}
  \geq 
  \frac{1}{k} 
 \sum_{\ell=1}^{k} V_{k\ell} 
  \Bigg]
 =
   h_{k,F}
,
\end{eqnarray*}
which establishes the result.   
  \cqfd
  \vspace{2mm}
 
 Note that this result shows that, for any~$F$ in~$\mathcal{F}$, the maximal depth,~$h_{k,F}$, in  Theorem~\ref{thetheorem} is monotone non-increasing in~$k$, 
 %. Since, irrespective of~$k$, the maximal depth is lower bounded by zero, the sequence of maximal depth 
 hence converges as~$k$ goes to infinity. Clearly, the law of large numbers and Slutzky's theorem imply that
% $$
%\frac{V_{k1}}{\frac{1}{k} \sum_{\ell=1}^k V_{k\ell}} 
%\to
%\frac{V_{11}}{{\rm E}[V_{11}]} 
%$$
$
V_{k1}/(\frac{1}{k} \sum_{\ell=1}^k V_{k\ell})
\to
V_{11}/{\rm E}[V_{11}]
$
in distribution, so that~$h_{k,F}$ converges to~$h_{\infty,F}=P[V_{11}\geq {\rm E}[V_{11}]]$ as~$k$ goes to infinity. In particular, for~$F=F_\alpha$, the limiting value is~$h_{\infty,F_\alpha}=P[Z_\alpha>\alpha]$, where~$Z_\alpha$ is~Gamma$(\alpha,1)$ distributed.  
\vspace{1mm}
 
We can now prove the main result of this paper. 
\vspace{1mm}

{\sc Proof of Theorem~\ref{thetheorem}}. 
We are looking for the infimum with respect to~$u=(u_1,\ldots,u_k)'$ in~$\mathcal{S}^{k-1}$, or equivalently in~$\R^k\setminus\{0\}$, of
\begin{eqnarray*}
\lefteqn{
p(u)
:=
P\Bigg[\sum_{\ell=1}^{k}u_\ell \Big( X_{k\ell}-\frac{1}{k}\Big) \geq 0 \Bigg]
=
P\Bigg[\sum_{\ell=1}^{k}u_\ell X_{k\ell}\geq \bar{u} \Bigg]
}
\\[2mm]
& & 
\hspace{3mm} 
=
P\Bigg[\sum_{\ell=1}^{k}u_\ell V_{k\ell}\geq \bar{u} \sum_{\ell=1}^k V_{k\ell} \Bigg]
=
P\Bigg[\sum_{\ell=1}^{k} \big(u_\ell-\bar{u}\big) V_{k\ell}\geq 0 \Bigg]
,	
\end{eqnarray*}
where we wrote~$\bar{u}:=\frac{1}{k}\sum_{\ell=1}^{k} u_\ell$. Without loss of generality, we may assume that~$u_1\geq u_2\geq\ldots\geq u_k$, which implies that~$u_1\geq \bar{u}$. Actually, if~$u_1=\bar{u}$, then all $u_\ell$'s must be equal to~$\bar{u}$, which makes the probability~$p(u)$ equal to one. Since this cannot be the infimum, we may assume that~$u_1>\bar{u}$, which implies that~$u_k< \bar{u}$. Therefore, denoting as~$m$ the largest integer for which~$u_m\geq \bar{u}$, we have~$1\leq m\leq k-1$.    
Then, letting~$s_m(u)=\sum_{\ell=1}^{m} \big(u_\ell-\bar{u}\big)$, we may then write
$$
p(u)
=
P\Bigg[
\sum_{\ell=1}^{m} \big(u_\ell-\bar{u}\big) V_{k\ell}
\geq 
\sum_{\ell=m+1}^{k} \big(\bar{u}-u_\ell\big) V_{k\ell}
 \Bigg]
=
P\Bigg[
\frac{
\sum_{\ell=m+1}^{k} d_\ell(u) V_{k\ell}
}
{
\sum_{\ell=1}^{m} c_\ell(u) V_{k\ell}
}
\leq 
1
 \Bigg]
,
$$
where
$
c_\ell(u)
=
(u_\ell-\bar{u})/s_m(u)
$,
$\ell=1,\ldots,m$
and
$
d_\ell(u)
=
(\bar{u}-u_\ell)/s_m(u)
$,
$\ell=m+1,\ldots,k$
are nonnegative and satisfy~$\sum_{\ell=1}^m c_\ell(u)=1$ and 
$$
\sum_{\ell=m+1}^k d_\ell(u)
=
\frac{
\sum_{\ell=m+1}^k (\bar{u}-u_\ell)
}
{
\sum_{\ell=1}^{m} \big(u_\ell-\bar{u}\big)
}
=
\frac{
(k-m) \bar{u} - \sum_{\ell=m+1}^k u_\ell
}
{
\sum_{\ell=1}^{m} u_\ell - m\bar{u}
}
= 1
.
$$
Since 
$\sum_{\ell=1}^{m} d_\ell(u) V_{k\ell}$ is unimodal at zero and since 
$
(c_{1}(u),c_{2}(u),\ldots,c_{m}(u))'$ is majorized by
$(1,0,\ldots,0)'\in\R^m$,
Lemma~\ref{thelemma} yields
$$
p(u)
\geq
P\Bigg[ 
\frac{
\sum_{\ell=m+1}^{k} d_\ell(u) V_{k\ell}
}
{
V_{k1}
}
\leq 
1
 \Bigg]
=
P\Bigg[ 
\frac{
V_{km}
}
{
\sum_{\ell=m+1}^{k} d_\ell(u) V_{k\ell}
}
\geq 
1
 \Bigg]
,
$$
where the lower bound is obtained for~$c_1(u)=1$ and~$c_2(u)=\ldots=c_m(u)=0$, that is, for~$u_1(>\bar{u})$ arbitrary and~$u_2=\ldots=u_m=\bar{u}$. Now, since 
$
(k-m)^{-1} 1_{k-m}$ is majorized by 
$(d_{m+1}(u),\ldots,d_{k}(u))'$ for any~$u$,
the same result provides 
\begin{eqnarray*}
\lefteqn{
\hspace{2mm} 
p(u)
\geq
P\Bigg[
\frac{
V_{km}
}
{
\frac{1}{k-m} \sum_{\ell=m+1}^{k} V_{k\ell}
}
\geq 
1
 \Bigg]
=
P\Bigg[
(k-m) V_{km}
\geq 
\sum_{\ell=m+1}^{k} 
V_{k\ell} 
 \Bigg]
}
\\[2mm]
& & 
\hspace{-5mm} 
=
P\Bigg[
(k-m+1) V_{km}
\geq 
\sum_{\ell=m}^{k} 
V_{k\ell} 
 \Bigg]
=
P\Bigg[
X_{k-m+1,1}
\geq  
\frac{1}{k-m+1} 
 \Bigg]
 = 
 h_{k-m+1,F}
,  
\end{eqnarray*}
with the lower bound obtained for~$d_{\ell}(u)=1/(k-m)$, $\ell=m+1,\ldots,k$, that is, for 
$
u_\ell = \bar{u}-\frac{1}{k-m} (u_1-\bar{u})
$, $\ell=m+1,\ldots,k$.
Therefore, the global minimum is the minimum of~$h_{k-m+1,F}$, $m=1,\ldots,k-1$,
which, in view of Lemma~\ref{thelemma}, is~$h_{k,F}$. This establishes the result. 
\cqfd
\vspace{3mm}

\begin{figure}
\begin{center}
	\includegraphics[width=.7\textwidth]{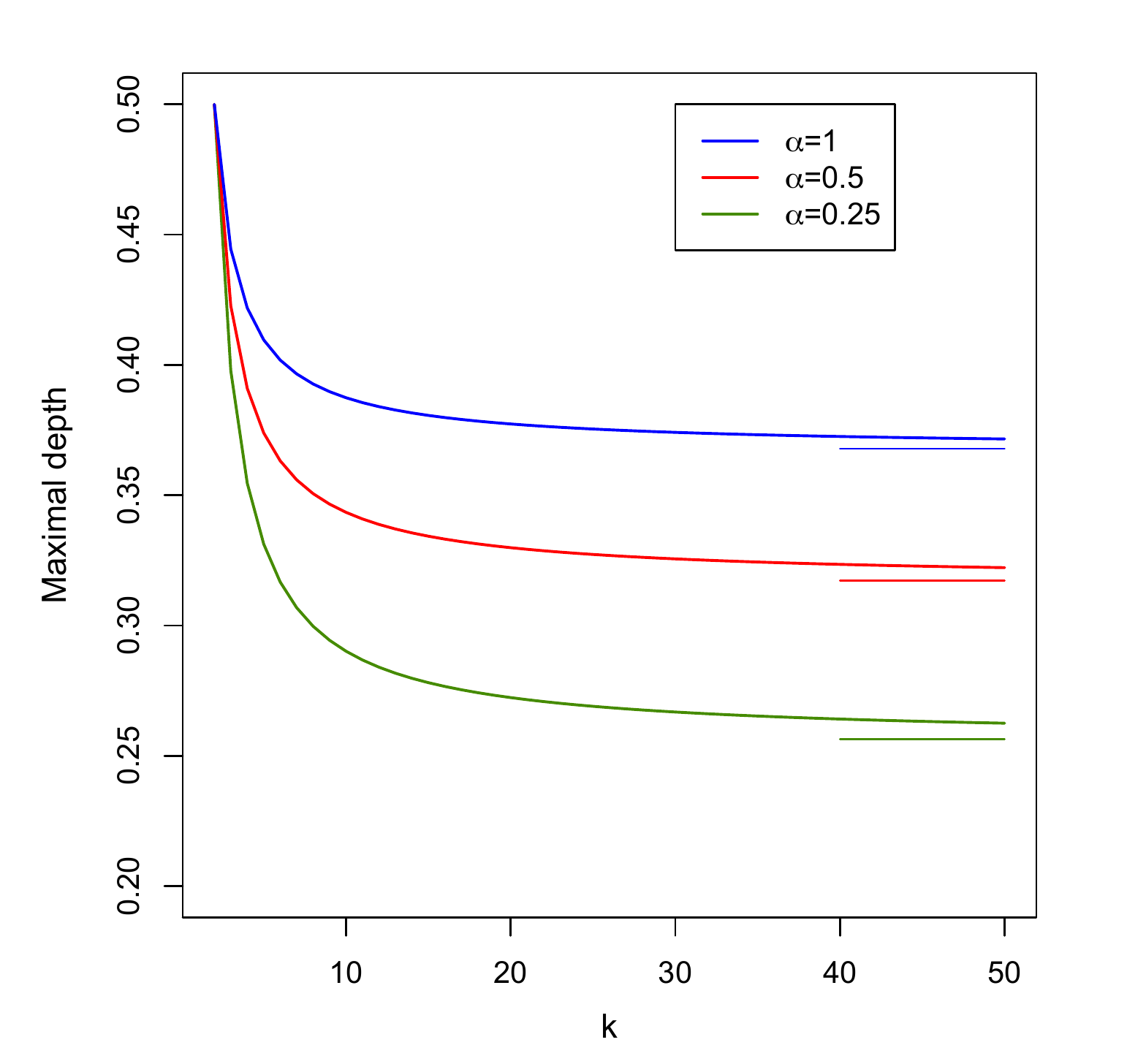}	
\end{center}
\caption{Plots of the maximal depth~$h_{k,F_\alpha}$ as a function of~$k$ for several values of~$\alpha$, where~$F_\alpha$ denotes the cumulative distribution function of the Gamma($\alpha,\frac{1}{2}$) distribution. For each value of~$\alpha$ considered, the limiting value as~$k$ goes to infinity is also showed.}
\label{Fig1}
\end{figure}

Figure~\ref{Fig1} plots the maximal depth~$h_{k,F_\alpha}$ as a function of~$k$ for several~$\alpha$, where~$F_\alpha$ still denotes the cumulative distribution function of the Gamma($\alpha,\frac{1}{2}$) distribution. In accordance with Lemma~\ref{thelemma}, the maximal depth is decreasing in~$k$ and is seen to converge to the limiting value~$h_{\infty,F_\alpha}$ that was obtained below that lemma. For any~$\alpha$, the maximal depth is equal to~$1/2$ if and only if~$k=2$, which is in line with the fact that the (non-atomic) probability measure~$P_{k,F_\alpha}$ is (angularly) symmetric about~$\mu_k$ if and only if~$k=2$; see \cite{RouStr2004}. Interestingly, thus, the asymmetry of~$P_{k,F}$ for~$k\geq 3$ would typically be missed by a test of symmetry that would reject the null when the sample (Tukey) median is too far from the sample mean.

%%%%%

\section{Numerical illustration}
\label{SecSimu}

\begin{figure}  
\begin{center}
	\includegraphics[width=.7\textwidth]{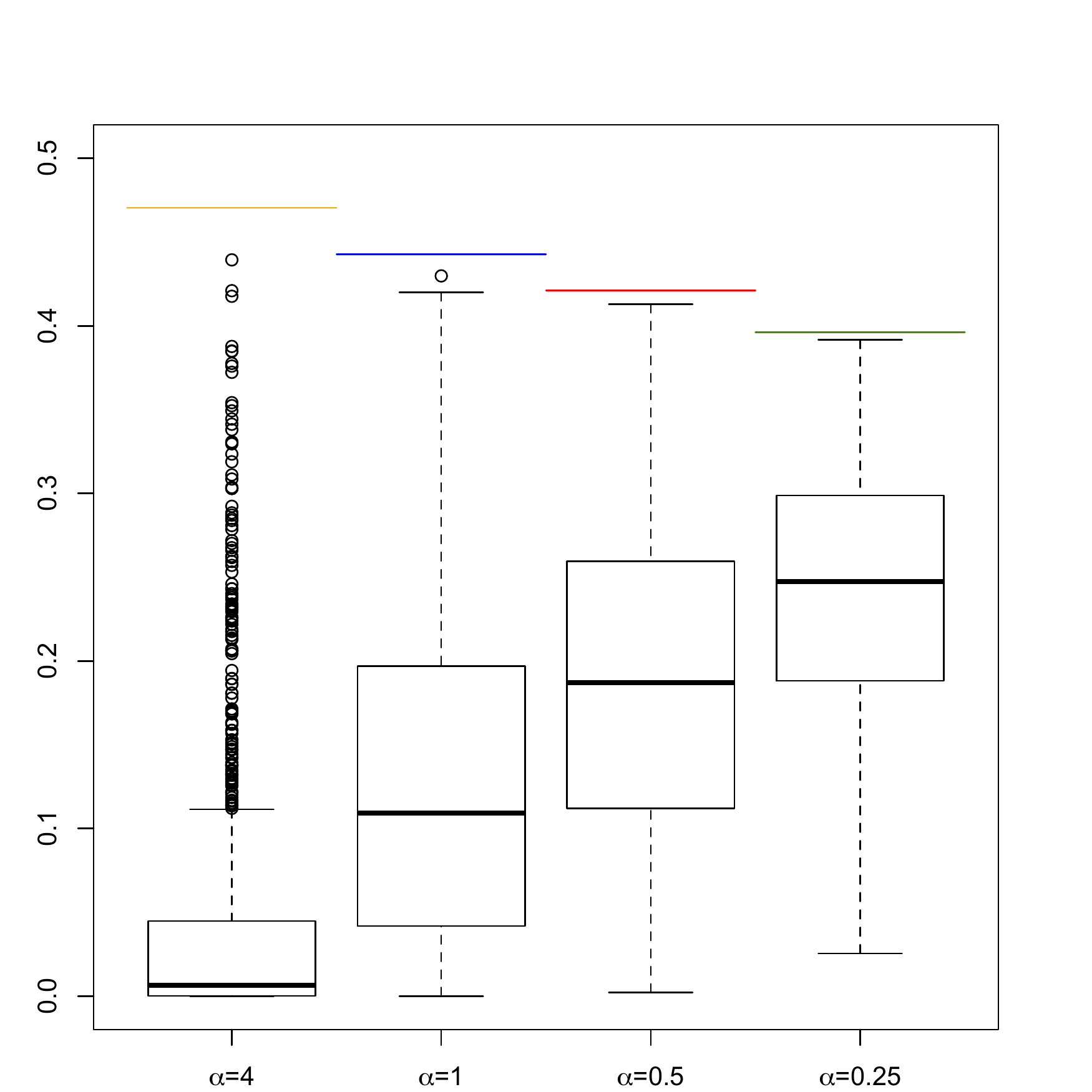}	 
\end{center}
\caption{Boxplots, for~$k=3$ and for various values of~$\alpha$, of the sample depth values~$H\!D(\theta_i,P_n)$, $i=1,\ldots,N$, of $N=1,000$ locations randomly drawn from the uniform distribution over the probability simplex~$\mathcal{D}_k$, where~$P_{n}$ denotes the empirical probability measure associated with a random sample of size~$n=100,000$ from the distributions~$P_{k,F_\alpha}$, with~$\alpha=4,1,0.5,0.25$. For each value of~$\alpha$ considered, the sample depth~$H\!D(\mu_k,P_n)$ is also provided (as explained in Section~\ref{SecSimu}, these sample depths were actually averaged over~$M=100$ mutually independent samples in each case).}
\label{Fig2}
\end{figure}

We conducted two numerical exercises to illustrate Theorems~\ref{thetheoremmax} and~\ref{thetheorem}, both in the trivariate case~$k=3$. In the first exercise, we generated~$N=1,000$ random locations~$\theta_1,\ldots,\theta_N$ 
%on~$\mathcal{D}_k$ 
from the uniform distribution over~$\mathcal{D}_k$ (that is the distribution~$P_{k,F_1}$ associated with the Gamma($1,\frac{1}{2}$) distribution of the~$V_{k\ell}$'s; see, e.g., Proposition~2 in \citealp{Bel2011}). Our goal is to compare, for various values of~$\alpha$, the depths~$H\!D(\theta_i,P_{k,F_\alpha})$, $i=1,\ldots,N$ with the 
%maximal 
depth~$H\!D(\mu_k,P_{k,F_\alpha})$.
%, where~$F_\alpha$ stands for the cumulative distribution function of the Gamma($\alpha,\frac{1}{2}$) distribution. 
For each~$\alpha$, these $N+1$ depth values were estimated by the depths~$H\!D(\theta_i,P_n)$, $i=1,\ldots,N$, and~$H\!D(\mu_k,P_n)$, computed with respect to the empirical measure~$P_n$ of a random sample of size~$n=100,000$ from~$P_{k,F_\alpha}$ (these~$N+1$ sample depth values were actually averaged over~$M=100$ mutually independent such samples). For each~$\alpha$, %this provides a collection of~$N+1$ depth values that are estimating~$H\!D(\theta_i,P_{k,F_\alpha})$, $i=1,\ldots,N$, and~$H\!D(\mu_k,P_{k,F_\alpha})$. 
Figure~\ref{Fig2} reports the boxplots of the resulting estimates of~$H\!D(\theta_i,P_{k,F_\alpha})$, $i=1,\ldots,N$, and marks the estimated value of~$H\!D(\mu_k,P_{k,F_\alpha})$. The results clearly support the claim in Theorem~\ref{thetheoremmax} that~$\mu_k=\arg\max_\theta H\!D(\theta,P_{k,F_\alpha})$.
% for any~$\alpha$.     

\begin{figure}
\begin{center}
	\includegraphics[width=1.01\textwidth]{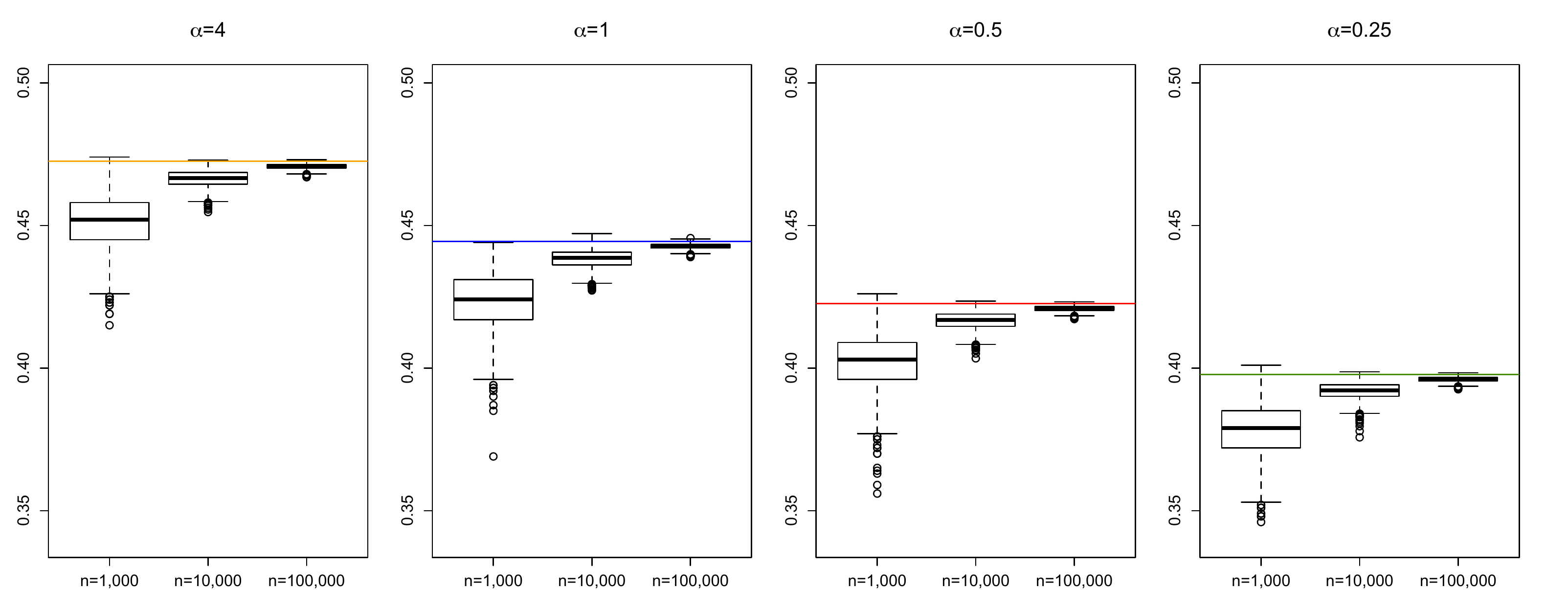}	
\end{center}
\caption{Boxplots, for~$k=3$ and for various values of~$\alpha$ and~$n$, of $M=1,000$ mutually independent values of the depth~$H\!D(\mu_k,P_{n})$, where~$P_{n}$ denotes the empirical probability measure associated with a random sample of size~$n$ from the same distributions~$P_{k,F_\alpha}$ as in Figure~\ref{Fig2}.}
\label{Fig3}
\end{figure}

In the second numerical exercise, we generated, for various values of~$\alpha$ and~$n$, a collection of $M=1,000$ mutually independent random samples of size~$n$ from the distribution~$P_{k,F_\alpha}$. For each sample, we evaluated the halfspace depth of~$\mu_k$ with respect to the corresponding empirical distribution~$P_n$. Figure~\ref{Fig3} provides, for each~$\alpha$ and~$n$, the boxplot of the resulting $M$ depth values. Clearly, the results, through the consistency of sample depth, support the theoretical depth values provided in Theorem~\ref{thetheorem}. 
Actually, while the theorem was only proved above for~$\alpha\leq 1$ (due to the unimodality condition in Lemma~\ref{Olshen}), these empirical results suggest that the theorem might hold also for~$\alpha>1$. The extension of the proof of Theorem~\ref{thetheorem} to further distributions is an interesting research question, that requires another approach (simulations indeed reveal that  Lemma~\ref{Olshen} does not hold if the cumulative distribution function of~$W$ is~$F_\alpha$ with~$\alpha>1$). Of course, another challenge is to derive a closed form expression for~$H\!D(\theta,P_{k,F})$ for an arbitrary~$\theta$. After putting some effort into this question, it seemed to us that such a computation calls for a more general stochastic ordering result than the one in Lemma~\ref{Olshen}.

%%%%%%%%%%%%%%%%%%%%%%%%%%%%%

\section*{Acknowledgements}

Davy Paindaveine's research is supported by the IAP research network grant \mbox{nr.} P7/06 of the Belgian government (Belgian Science Policy), the Cr\'{e}dit de Recherche  J.0113.16 of the FNRS (Fonds National pour la Recherche Scientifique), Communaut\'{e} Fran\c{c}aise de Belgique, and a grant from the National Bank of Belgium. Germain Van Bever's research is supported by the FC84444 grant of the FNRS. % (Fonds National pour la Recherche Scientifique), Communaut\'{e} Fran\c{c}aise de Belgique.

%%%%%%%%%%%%%%%%%%%%%%%%%%%%%%%%%%%

%\newpage
\bibliographystyle{elsarticle-harv}
\bibliography{Paper}

\end{document}